\documentclass[10pt]{article}
\usepackage{amsmath}
\usepackage{latexsym}
\usepackage{amssymb}
\usepackage{amsfonts}
\usepackage{amsthm}
\title{COMPUTER ALGEBRA AND LANCZOS POTENTIAL}
\author{J.-F. Pommaret \\ CERMICS, Ecole des Ponts ParisTech, France \\
 jean-francois.pommaret@wanadoo.fr \\
 (http://cermics.enpc.fr/$\sim$pommaret/home.html )}
\date{  }
\begin{document}
\maketitle

\noindent
{\bf ABSTRACT}  \\
 We found in 2016 a few results on the mathematical structure of the conformal Killing differential sequence in arbitrary dimension $n$, in particular the rank and order changes of the successive differential operators for $n=3,n=4$ or $n\geq 5$. They were so striking that we did not dare to publish them before our former PhD student A. Quadrat (INRIA) could confirm them while using new computer algebra packages that he developped for studying extension modules in differential homological algebra. In the meantime, as a complementary result, we found in 2017 the "{\it missing link} " justifying the doubts we had since a long time on the origin and existence of Gravitational Waves in General Relativity. In both cases, the main tool is the explicit computation of certain extension modules for the classical or conformal Killing differential sequences. These results therefore lead to revisit the work of C. Lanczos and successors on the existence of a parametrization of the Riemann or Weyl operators and their  respective formal adjoint operators. We also provide an example showing how these extension modules are depending on the {\it structure constants} appearing in the {\it Vessiot structure equations } (1903), still not acknowledged after one century even though they generalize the constant Riemannian curvature integrability condition of L.P. Eisenhart (1926) for the Killing equations. The present paper is written from a lecture gven at the recent 24 th conference on Applications of Computer Algebra (ACA 2018) held in Santiago de Compostela, Spain, june 18-22, 2018.  \\  \\

\vspace{3cm}


\newpage

\noindent
{\bf 1) INTRODUCTION}  \\

We start this paper with a brief synthetic review introducing a few notations allowing to describe the standard procedure used in any book on General relativity (GR) in order to justify the existence of gravitational waves. A main idea, rarely pointed out, is the concept of {\it linearization} ("{\it lin}" for short) and its systematic use in mathematical physics. When $X$ is a manifold of dimension $n$ and ${\cal{E}}$ is a fibered manifold over $X$, we may introduce the {\it vertical bundle} $V({\cal{E}})$ as a vector bundle over ${\cal{E}}$ and, for any section $f$ of ${\cal{E}}$, introduce the reciprocal image $f^{-1}(V({\cal{E}}))$ as a vector bundle over $X$ with the same fiber dimension ([13],[15],[21]). In the formulas below, ${\cal{L}}$ is the {\it Lie derivative} and the order of a differential operator is written under its representative arrow. All the results presented are formal and local but the corresponding global notations will be used for simplicity. \\
 
\noindent
\fbox{1} \,\,\, {\bf METRIC}:  \,\,\,$\omega=({\omega}_{ij})=({\omega}_{ji})\in S_2T^*, det(\omega)\neq 0 \,\,\, \stackrel{lin}{\longrightarrow}\,\,\, \Omega \in {\omega}^{-1}(V(S_2T^*))=S_2T^*$\\  \\

\noindent
\fbox{2} \,\,\,  {\bf CHRISTOFFEL SYMBOLS}: \,\,\,${\gamma}^k_{ij}=\frac{1}{2}{\omega}^{kr}({\partial}_i{\omega}_{rj} +{\partial}_j{\omega}_{ir} -{\partial}_r{\omega}_{ij}) \,\,\, \stackrel{lin}{\longrightarrow}  \,\,\, \Gamma \in S_2T^*\otimes T$  \\  \\

\noindent
\fbox{3} \,\,\, {\bf RIEMANN TENSOR}:\,\,\, ${\rho}^k_{l,ij}={\partial}_i{\gamma}^k_{lj} - {\partial}_j{\gamma}^k_{li}+{\gamma}^r_{li}{\gamma}^k_{rj} - {\gamma}^r_{lj}{\gamma}^k_{ri} \,\,\, \stackrel{lin}{\longrightarrow}  \,\,\, R \in {\wedge}^2T^*\otimes T^*\otimes T $ \\

\[  {\omega}_{kr}{\rho}^r_{l,ij}={\rho}_{kl,ij}= - {\rho}_{lk,ij}= - {\rho}_{kl,ji}={\rho}_{ij,kl} \Rightarrow {\rho}^r_{r,ij}=0 \]
\[    {\rho}^k_{l,ij}+{\rho}^k_{i,jl} + {\rho}^k_{j,li}=0  \Rightarrow {\rho}^r_{i,rj}={\rho}^r_{j,ri} \]   \\

\noindent
\fbox{4} \,\,\, {\bf KILLING OPERATOR}: \,\,\, ${\cal{D}}:T \rightarrow S_2T^*: \xi \rightarrow {\cal{L}}(\xi)\omega=\Omega  $\\

${\cal{D}}\xi=0, {\cal{D}}\eta=0 \Rightarrow {\cal{D}}[\xi,\eta]=0 $ \,\,\, ({\it Lie} operator)   \\

\noindent
The corresponding first order system $R_1=R_1(\omega)\subset J_1(T)$ with {\it symbol} $g_1=R_1\cap T^*\otimes T \subset J_1(T)$ is defined by the equations (See [13,15,28] for more details on jet theory):  \\

\noindent
\fbox{5} \,\,\, {\bf KILLING EQUATIONS}: \,\,\, $({\cal{L}}(\xi)\omega)_{ij}\equiv{\Omega}_{ij} \equiv {\omega}_{rj}{\partial}_i{\xi}^r +{\omega}_{ir}{\partial}_j {\xi}^r  +  {\xi}^r{\partial}_r{\omega}_{ij}=0  $ \\

\noindent
\fbox{6} \,\,\, {\bf KILLING SYMBOL}: \,\,\,  ${\omega}_{rj}v^r_i +{\omega}_{ir}v^r_j=0  $  \\ \\
\noindent
Looking for successive {\it compatibility conditions} (CC) we may exhibit:  \\

\noindent
\fbox{7} \,\,\, {\bf KILLING SEQUENCE}: \,\,\, $\Theta $ = \{ {\it Killing} vector fields \} \\
\[  \begin{array}{ccccccccccc}
 0 &\rightarrow &\Theta &\rightarrow &T & \underset 1{\stackrel{{\cal{D}}}{\longrightarrow}} &S_2T^* &\underset 2{\stackrel{{\cal{D}}_1}{\longrightarrow}}  &  F_1 & \underset 1{\stackrel{{\cal{D}}_2}{\longrightarrow }}& F_2    \\
     & & & & & & & & & &   \\
  & && &n & {\stackrel{Killing}{\longrightarrow}} &\frac{n(n+1)}{2} &{\stackrel{Riemann}{\longrightarrow}}  &  \frac{n^2(n^2-1)}{12} & {\stackrel{Bianchi}{\longrightarrow }}& \frac{n^2(n^2-1)(n-2)}{24} 
\end{array}    \]

\noindent
\fbox{8} \,\,\, {\bf RICCI TENSOR}: \,\,\,  ${\rho}_{ij}={\rho}^r_{i,rj}={\rho}_{ji}  \,\,\, \stackrel{lin}{\longrightarrow}  \,\,\, (R_{ij})\in S_2T^*$ \\  \\

\noindent
\fbox{9} \,\,\,  {\bf EINSTEIN TENSOR}: \,\,\,  \fbox{${\epsilon}_{ij}={\rho}_{ij} - \frac{1}{2}{\omega}_{ij} {\omega}^{rs}{\rho}_{rs}$} \,\,\,  $ \stackrel{lin}{\longrightarrow}\,\,\, (E_{ij})\in S_2T^*$  \\  \\

\noindent
\fbox{10} \,\,\,  {\bf EINSTEIN EQUATIONS}:  \,\,\, $(Einstein) \,\,\,div(E)=0 \Rightarrow  $ \fbox{ $E_{ij} \sim {\Sigma}_{ij}$}$\Leftarrow 
div({\Sigma})= 0  \,\,\,(Cauchy) $   \\  \\

\noindent
\fbox{11} \,\,\, {\bf WAVE EQUATIONS}: \,\,\, \fbox{${\bar{\Omega}}_{ij}={\Omega}_{ij} - \frac{1}{2} {\omega}_{ij}{\omega}^{rs}{\Omega}_{rs}$}  \,\,\, $\Rightarrow \,\,\, $ \fbox{ $ \Box {\bar{\Omega}}_{ij}+ ... \sim {\Sigma}_{ij}$}\\  

As we shall see, these results, which are of course mathematically correct, are not conceptually coherent with differential homological algebra and the well known Poincar\'{e} duality existing between GEOMETRY (Killing sequence) and PHYSICS ( adjoint sequence). Such results, based on the search for parametrizing an operator, will AUTOMATICALLY lead to revisit the work of C. Lanczos ([8],[9]) for parametrizing {\it Riemann} or {\it Weyl} operators both with their formal adjoint operators, explaining in particular why the many papers written by his followers during more than $50$ years are containing so many contradictory claims ([1-4],[10],[12]).  \\ \\

\noindent 
{\bf 2) PARAMETRIZATION}  \\

Looking for a  differential sequence of the form  \,\,\, $ \xi  \stackrel{{\cal{D}}}{\longrightarrow } \eta \stackrel{{\cal{D}}_1}{\longrightarrow} \zeta  $  
\,\,\, in which ${\cal{D}}_1$ generates all the {\it compatibility conditions} (CC) of ${\cal{D}}$ and $(\xi,\eta,\zeta)$ are sections of certain vector bundles, two problems can be proposed at once:  \\  \\

\noindent
\fbox{ {\bf DIRECT PROBLEM}} \,\,\,\,  ${\cal{D}}$ given, find ${\cal{D}}_1$ : M. Janet (1920), D.C. Spencer (1970)   \\

\hspace*{1cm} $\Updownarrow  $  \\

\noindent
{\fbox { {\bf INVERSE PROBLEM}} \,\, $ {\cal{D}}_1$ given, find ${\cal{D}}$ : not always possible  \\  

Before solving negatively in ([14]) the challenge of J. Wheeler (1970) concerning the possibility to parametrize Einstein equations (See also [30]), we quote a few references on the direct problem ([ 6],[13],[15],[16],[21],[24],[28]) or the inverse problem ([16],[17],[30]) and present below the best elementary but absolutely non-trivial example we know from control theory: Ê\\ 

\noindent
{\bf THEOREM 2.1}: When $n=1$, a classical control system is controllable if and only if it is {\it parametrizable}, that is if and only if it generates the CC  of a previous operator describing therefore {\it a} parametrization} as there may be many different ones.  \\  

\noindent
 {\bf COROLLARY 2.2}: A classical control system defined over an ordinary differential field $K$ by equations linearly independent over the ring 
 $D=K[d]$ of differential operators with coefficients in $K$ is controllable if and only if the formal adjoint of the corresponding differential operator is injective, even when $D$ is non-commutative.  \\ 

\noindent
{\bf EXAMPLE 2.3}: {\it DOUBLE PENDULUM}\\

Rigid bar of length $L$ moving along the left to right horizontal axis $0x$, downwards vertial axis $0y$ parallel to gravity $g$, first pendulum made by a mass $m_1$, having length $l_1$ and moving by an angle ${\theta}_1$ with respect to the vertical, second pendulum made by a mass $m_2$, having length $l_2$ and mouving by an angle ${\theta}_2$ with respect to the vertical.\\
\noindent

Control system:  \fbox{ $\ddot{x} + l_1 {\ddot{\theta}}_1 + g {\theta}_1=0, \,\,\, \,\,\ddot{x} + l_2 {\ddot{\theta}}_2 + g {\theta}_2=0 $ }\\  

\noindent
{\it Parametrization}:\\
\noindent
$\bullet$ $l_1\neq l_2$:  \\
\[ \left\{  \begin{array}{rcl}
  - l_1 l_2 d^4\phi - g(l_1+l_2)d^2 \phi - g^2 \phi & = & x  \\
  l_2 d^4\phi +g d^2 \phi & = & {\theta}_1  \\
  l_1 d^4 \phi + gd^2 \phi  & = & {\theta}_2
  \end{array} \right.  \]  \\
  
 \noindent 
 $\bullet$ $l_1=l_2 =l, \theta={\theta}_1 - {\theta}_2 \Rightarrow l\ddot{\theta} + g \theta=0$, \,\, $\theta(0)=0, \dot{\theta}(0)=0 \Rightarrow \theta (t)=0$.  \\  
 
Multiplying on the left the two OD equations by two {\it test functions} $({\lambda}^1,{\lambda}^2)$ and integrating by parts, we let the reader check by himself the second theorem as follows:   \\
\[ {\ddot{\lambda}}^1 +{\ddot{\lambda}}^2=0,\,\,  l_1 {\ddot{\lambda}}^1+ g {\lambda}^1=0, \,\, l_2 {\ddot{\lambda}}^2 + g{\lambda}^2=0 \,\, \stackrel{l_1\neq l_2}{\Longleftrightarrow} \,\, {\lambda}^1=0,\,\, {\lambda}^2=0  \]  
We finally notice that  COMPUTER ALGEBRA IS ABSOLUTELY NEEDED for treating the case $l_1=cst, l_2=l_2(t)$ ([16]).  \\  \\

\noindent                                                                                                                                      
{\bf  3) DOUBLE DUALITY TEST}  \\

The test is based on a systematic use of the (formal) {\it adjoint} ($ad$) of an operator and has five steps described in the diagram below ([7,15,16,24]): \\

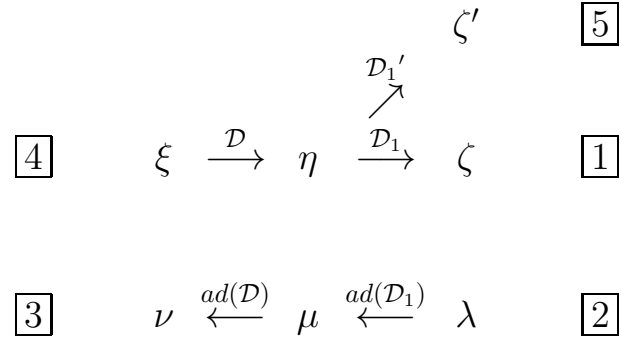
  \[  \Large {\begin{array}{rcccccl}
 & & & & &  {\zeta}' &\hspace{1cm} \fbox{5}  \\
 & & & & \stackrel{{{\cal{D}}_1}'}{\nearrow} &  &  \\
\fbox{4} \hspace{1cm}& \xi  & \stackrel{{\cal{D}}}{\longrightarrow} &  \eta & \stackrel{{\cal{D}}_1}{\longrightarrow} & \zeta &\hspace{1cm}   \fbox{1}  \\
 &  &  &  &  &  &  \\
 &  &  &  &  &  &  \\
 \fbox{3} \hspace{1cm}& \nu & \stackrel{ad({\cal{D}})}{\longleftarrow} & \mu & \stackrel{ad({\cal{D}}_1)}{\longleftarrow} & \lambda &\hspace{1cm} \fbox{2}
  \end{array}  }  \] 
\vspace{2mm}  \\

  \[ \large{\begin{tabular}{c}
  $ ad(ad({\cal{D}}))={\cal{D}}, \,\,\,  ad({\cal{D}})\circ ad({\cal{D}}_1)= ad ({\cal{D}}_1 \circ {\cal{D}})=0   $ \\  \\
   $\Rightarrow   {\cal{D}}_1\circ {\cal{D}} =0 \Rightarrow {\cal{D}}_1 \,\, AMONG \,\, the \,\,CC \,\, of \,\, {\cal{D}} $ \\  \\
Step \fbox{5} $\Rightarrow {{\cal{D}}_1}' \,\, GENERATES \,\, the \,\, CC \,\, of \,\, {\cal{D}} \,\,\,  \Rightarrow \,\,\, $ \fbox{ ${\cal{D}}_1 \leq 
{{\cal{D}}_1}' $ } 
\end{tabular}  }  \] \\

\noindent
{\bf THEOREM 3.1}: ${\cal{D}}_1$ parametrized by ${\cal{D}}$ \,\, $ \Leftrightarrow $ \,\, \fbox{${\cal{D}}_1={{\cal{D}}_1}'$}. \\ 
 
\noindent
{\bf COUNTEREXAMPLE 3.2}:  \fbox{EINSTEIN EQUATIONS CANNOT BE PARAMETRIZED }  \\
Contrary to the Ricci operator ($4$ terms only), the Einstein operator ($6$ terms) is SELF-ADJOINT, the sixth terms being exchanged between themselves under $ad$:  
\[ \Huge { {\lambda}^{ij}({\omega}^{rs}d_{ij}{\Omega}_{rs}) \stackrel{ad}{\longleftrightarrow} ({\omega}^{rs}d_{ij}{\lambda}^{ij}){\Omega}_{rs}=({\omega}_{ij}d_{rs}{\lambda}^{rs}){\Omega}^{ij} } \] \\ 

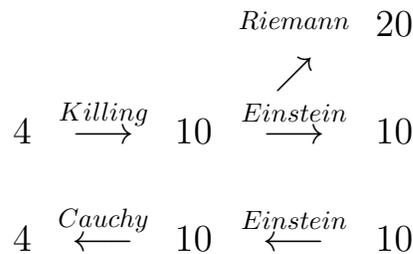
\[ \Large { \begin{array}{rcccl}
  &  &  &\stackrel{Riemann}{ }  & 20   \\
  & &  & \nearrow &    \\
 4 &  \stackrel{Killing}{\longrightarrow} & 10 & \stackrel{Einstein}{\longrightarrow} & 10  \\
  & & & &  \\
 4 & \stackrel{Cauchy}{\longleftarrow} & 10 & \stackrel{Einstein}{\longleftarrow} &  10  
\end{array} } \]

      \newpage

\noindent
{\bf 4) VARIATIONAL CALCULUS WITH CONSTRAINTS }    \\  

\noindent
In order to explain the motivation of Lanczos in ([8]) for using the {\it Lagrange multiplier} $\lambda$, let us suppose that ${\cal{D}}_1$ generates the CC of  ${\cal{D}}$ {\it AND} that $ad({\cal{D}})$ generates the CC of $ad({\cal{D}}_1)$. There are two possible points of view, exactly like in continuum mechanics:   \\

\[ \begin{array}{rcl}
\fbox{${\cal{D}}\xi=\eta $ }\,\, \Rightarrow 
\hspace{25mm}  \Phi={\int}_V\varphi (\eta)dx  \Rightarrow   \delta \Phi  &= &  \int \frac{\partial \varphi}{\partial \eta}\delta \eta dx          \\
                                                                                  & = &  \int \frac{\partial \varphi}{\partial \eta}{\cal{D}}\delta \xi  dx  \\
                                                                                  &  =  & \int (ad({\cal{D}})\frac{\partial \varphi}{\partial \eta})\delta \xi dx +... \\
                                                                                   & \Rightarrow &  \,\,\, Cauchy
\end{array}   \]

\[  \begin{array}{rcl}
\fbox{${\cal{D}}_1\eta=0 $}\,\, \Rightarrow 
\hspace{5mm}  \Phi={\int}_V(\varphi (\eta) - \lambda {\cal{D}}_1\eta)dx  \Rightarrow   \delta \Phi  &= &  \int (\frac{\partial \varphi}{\partial \eta}\delta \eta - \lambda {\cal{D}}_1\delta \eta)dx          \\
                                                                                  & = &  \int (\frac{\partial \varphi}{\partial \eta}- ad({\cal{D}}_1)\lambda)\delta \eta  dx + ...
                                                                                 
\end{array}   \]    

\[  \Rightarrow  \mu=\frac{\partial \varphi}{\partial \eta}=ad({\cal{D}}_1)\lambda \,\, (parametrization \,\, by \,\,  \lambda) \stackrel{ad({\cal{D}})}{\longrightarrow} ad({\cal{D}})\mu=0  \,\, (elimination \,\, of \,\, \lambda)\]  \\

\noindent
{\bf EXAMPLE 4.1}: \,\, {\fbox { $n=2$}} \,\,\,  Airy parametrization (1863)

\[  \begin{array}{ccccccl}
2 & \stackrel{Killing}{\longrightarrow}& 3 & \stackrel{Riemann}{\longrightarrow} & 1 & \rightarrow & 0  \\
  &  &  &  &  &  &  \\
2 & \stackrel{Cauchy}{\longleftarrow} & 3 & \stackrel{Airy}{\longleftarrow} & 1 &  &
\end{array}   \]   \\

\noindent
{\bf THEOREM 4.2}:  \hspace{1cm}  $Cauchy=ad(Killing)$, \hspace{5mm} $Airy=ad(Riemann)$

\[  \lambda (d_{22}{\Omega}_{11} - \fbox{2} d_{12}{\Omega}_{12} + d_{11}{\Omega}_{22})=(d_{22}\lambda{\Omega}_{11} - \fbox{2} d_{12}\lambda{\Omega}_{12}+d_{11}\lambda{\Omega}_{22}) + ...  \]

\[ {\sigma}^{ij}{\Omega}_{ij}= {\sigma}^{11}{\Omega}_{11} + \fbox{2} {\sigma}^{12}{\Omega}_{12}+ {\sigma}^{22}{\Omega}_{22}  \]

\noindent
  Cauchy \,\,\, \fbox{$d_1{\sigma}^{11} + d_2{\sigma}^{12}=f^1, \,\,\, d_1{\sigma}^{21} + d_2{\sigma}^{22}=f^2$}  \\

\noindent
Airy \,\,\, \fbox{ ${\sigma}^{11}=d_{22}\lambda, \,\,\,{\sigma}^{12}={\sigma}^{21}= - d_{12} \lambda, \,\,\, {\sigma}^{22}=d_{11}\lambda $}\\

It is an open problem to know why one may sometimes find a SELF-ADJOINT OPERATOR:  \\

\noindent
{\bf EXAMPLE 4.3}: \,\, \fbox{n=3}   \,\,\,\,We now present the  BELTRAMI PARAMETRIZATION (1892):  \\

\noindent
  \[ \large { \left\{  \begin{array}{c}
{\sigma}^{11}  \\  {\sigma}^{12} \\ {\sigma}^{12} \\ {\sigma}^{13} \\ {\sigma}^{23} \\  {\sigma}^{33}
\end{array}  \right\} =   \left\{  \begin{array}{cccccc}
0 & 0 & 0 & d_{33} & -2d_{23} & d_{22}  \\
0 & -d_{33} & d_{23} & 0 & d_{13} & - d_{12}  \\
0 & d_{23} & - d_{22} & - d_{13} & d_{12} & 0  \\
d_{33} & 0 & - 2 d_{13} & 0  &  0  & d_{11}   \\
- d_{23} & d_{13} & d_{12} & 0 &  d_{11} & 0  \\
d_{22} &- 2d_{12}& 0 & d_{11} &  0 & 0 
\end{array}  \right\}  \left\{\begin{array}{c} {\phi}_{11}\\{\phi}_{12}\\ {\phi}_{13}\\ {\phi}_{22} \\ {\phi}_{23} \\ {\phi}_{33}
\end{array}  \right\} } \]  \\ 

\noindent
which  does not seem to be SELF-ADJOINT . \\

 \[ \large  {     \Leftrightarrow  \,\,\,d_r {\sigma}^{ir}=0   \,\,\,  (Cauchy)  } \]   
  
\noindent
Accordingly, the Beltrami parametrization of the {\it Cauchy} operator for the stress is nothing else than the formal adjoint of the {\it Riemann} operator, namely:  \\

\hspace*{4cm} \fbox{$ad(Riemann)=Beltrami$}  \\  \\

However, modifying slightly the rows, we get the new operator matrix:  \\

\noindent
\[ \large { \left\{  \begin{array}{c}
{\sigma}^{11}  \\  2{\sigma}^{12} \\ 2 {\sigma}^{13} \\ {\sigma}^{22} \\ 2{\sigma}^{23} \\  {\sigma}^{33}
\end{array}  \right\} =   \left\{  \begin{array}{cccccc}
0 & 0 & 0 & d_{33} & -2d_{23} & d_{22}  \\
0 & - 2 d_{33} & 2d_{23} & 0 & 2d_{13} & - 2d_{12}  \\
0 & 2d_{23} & - 2d_{22} & - 2d_{13} & 2d_{12} & 0  \\
d_{33} & 0 & - 2 d_{13} & 0  &  0  & d_{11}   \\
- 2d_{23} & 2d_{13} & 2d_{12} & 0 &  d_{11} & 0  \\
d_{22} &- 2d_{12}& 0 & d_{11}&  0 & 0 
\end{array}  \right\}  \left\{\begin{array}{c} {\phi}_{11}\\{\phi}_{12}\\ {\phi}_{13}\\ {\phi}_{22} \\ {\phi}_{23} \\ {\phi}_{33}
\end{array}  \right\} } \]     \\

\noindent
which is indeed SELF-ADJOINT .  \\

We end this section by noticing that MOST TEXTBOOKS, using the infinitesimal deformation tensor  $(\frac{1}{2}{\Omega}_{ij})$ in the Helmholtz free energy, are claiming that ${\sigma}^{ij}= 2 \frac{\partial \varphi}{\partial {\Omega}_{ij}}  \Rightarrow  {\sigma}^{ij}={\sigma}^{ji}$ but, as we have seen, such a claim is not correct at all (See [21] for the problems brought while using computer algebra computer with these tricky factors "$2$" involved).  \\\\

\noindent
{\bf 5) CONTRADICTIONS}  \\

Coming back to the mathematical origin of gravitational waves, we obtain:  \\

  \hspace*{3cm}    \fbox{   \large     {\bf FIRST CONTRADICTION} } \\

\noindent
\fbox {n=4}  \hspace {2cm}  $ ad(Killing)=Cauchy, \,\,\,  ad(Riemann)=Beltrami $ \\   

\[   \begin{array}{rcccccccccl}
 & \fbox{4} & \stackrel{Killing}{\longrightarrow} & 10 & \stackrel{Riemann}{\longrightarrow} & 20 & \stackrel{Bianchi}{\longrightarrow} & 20 & \longrightarrow & 6 & \rightarrow 0 \\
  &   &                                                            & \parallel &  \searrow     & \downarrow  &  & \downarrow &  &   \\
 &    &                                                            &  10     & \stackrel{Einstein}{\longrightarrow}  & 10  & \stackrel{div}{\longrightarrow} & \fbox{4 }   & \longrightarrow & 0                      &   \\  
 & & &  &   &   \downarrow  & &  \downarrow  &  &  &    \\                                                                                                                                                                                                                                                                                                                                                                                                                                                                                                                                                                                                                                                                                                                                                                                                                                                                                                                                                                                                                                                                                                                                                                                                                                                                                                     
  &   &  &  & & 0 &  & 0 &  &  \\   \\

0 \leftarrow & \fbox{4} & \stackrel{Cauchy}{\longleftarrow} & 10 & \stackrel{Beltrami}{\longleftarrow} & 20 & \stackrel{ad(Bianchi)}{\longleftarrow} & 20 &  & & \\
                     &           &                                                         & \parallel &  \nwarrow    & \uparrow  &  & \uparrow  &  &  \\
  &    &      & 10 &  \stackrel{Einstein}{\longleftarrow} & 10 & \stackrel{ad(div)}{\longleftarrow} &  \fbox{4} &   &  &  \\
  & & & & & \uparrow & &\uparrow  & & &   \\
  
  &  &  &  &  & 0 & & 0&  
\end{array}   \]  
Comparing the two above diagrams, it becomes  clear that the {\it Cauchy} operator has NOTHING TO DO with the {\it div} operator classically induced from the {\it Bianchi} operator by contracting indices because the \fbox{4} on the left has has a quite different mathematical meaning than the \fbox{4} on the right.  \\

  \hspace*{3cm}      \fbox{   \large     {\bf SECOND CONTRADICTION}}  

\[  \large {  E_{ij}=R_{ij}- \frac{1}{2}{\omega}_{ij}{\omega}^{rs}R_{rs} \,\,\, \Rightarrow \,\,\, E=C \circ R  }\]
\[  \large {  E: \Omega  \stackrel{C}{\longrightarrow}  \bar{\Omega}=\Omega - \frac{1}{2}\omega tr(\Omega) \stackrel{X}{\longrightarrow} 
S_2T^* } \] 

 $Einstein$ operator  $E$  (6 terms) $\rightarrow$ $wave$ operator $ X $  (4 terms only) 
\[  \large {  E=X\circ C \Rightarrow E=ad(E)=ad(C)\circ ad(X)=C \circ ad(X)=C \circ R } \] 
\[  \large { \Rightarrow ad(X)=Ricci \,\, \Rightarrow \,\,  \fbox{X=ad(Ricci)} } \]

\noindent
We may therefore surprisingly conclude by saying that:  \\

\noindent
THE EINSTEIN OPERATOR IS USELESS while ONLY THE RICCI OPERATOR IS USEFUL.  \\  \\

\noindent
{\bf 6) LANCZOS POTENTIAL}   \\

\noindent
We first notice that the {\it Poincar\'{e} sequence} for the exterior derivative $d$, namely: \\
\[  {\wedge}^0 T^* \stackrel{d}{\longrightarrow} {\wedge}^1T^* \stackrel{d}{\longrightarrow}{\wedge}^2T^* \stackrel{d}{\longrightarrow}  ....
\stackrel{d}{\longrightarrow} {\wedge}^nT^* \rightarrow 0  \]
IS SELF ADJOINT UP TO SIGN : \\

\noindent
{\bf EXAMPLE 6.1}: \fbox {n=3} \,\,\,\,  ${\wedge}^0 T^* \stackrel{grad}{\longrightarrow} {\wedge}^1T^* \stackrel{curl}{\longrightarrow}{\wedge}^2T^* \stackrel{div}{\longrightarrow} {\wedge}^3T^* \rightarrow 0  $
\[  ad(grad)=- div, \,\, ad(curl)=curl, \,\,\, ad(div)=- grad  \]

However, if we start from the following sequence and its formal ajoint:  \\

\[  \begin{array}{ccccc}
\xi &\longrightarrow & \left\{ \begin{array}{rcl} d_{22} \xi & = & {\eta}^2 \\  d_{12}\xi & = & {\eta}^1 \end{array} \right. & \longrightarrow &d_1{\eta}^2-d_2{\eta}^1= \zeta  \\
\xi & \stackrel{{\cal{D}}}{\longrightarrow} & \eta & \stackrel{{\cal{D}}_1}{\longrightarrow} & \zeta  \\
  &  &  &  &   \\
  \nu & \stackrel{ad({\cal{D}})}{\longleftarrow} & \mu & \stackrel{ad({\cal{D}}_1)}{\longleftarrow } & \lambda  \\
   d_{12}{\mu}^1 + d_{22}{\mu}^2=\nu                                        & \longleftarrow & \left\{  \begin{array}{rl} -d_1\lambda & = {\mu}^2 \\ d_2\lambda & ={\mu}^1 \end{array} \right. & \longleftarrow & \lambda  \\
            & \swarrow & & &     \\
            d_1{\mu}^1 +d_2{\mu}^2={\nu}' &   & & &                              
\end{array}  \]   \\
then, even if ${\cal{D}}_1$ generates the CC of ${\cal{D}}$, in general $ad({\cal{D}})$ may not generate all the CC of $ad({\cal{D}}_1)$. Such a "{\it GAP} ", namely the lack of formal exactness of the adjoint sequence when the initial sequence is formally exact, led to introduce the {\it extension modules} ([11],[27]) because of the following (difficult) theorems (See [11,27] or [22,23,24] for more details):  \\

\noindent
{\bf THEOREM 6.2}: If $M$ is the differential module defined by ${\cal{D}}$, the {\it extension modules} $ext^i(M)$ do not depend on the sequence used for their computation.  \\

\noindent
{\bf THEOREM 6.3}: The {\it Spencer sequence} for any Lie operator ${\cal{D}}$ which is coming from a Lie group of transformations is (locally) isomorphic to the tensor product of the Poincar\'{e} sequence by the corresponding finite Lie algebra ${\cal{G}}$.  \\

\noindent
{\bf COROLLARY 6.4}: \,\,\,  $ {ext}^1(M)=0,\,\, {ext}^2(M)=0 , ...  \Rightarrow $ \fbox{ NO GAP }\\

\noindent
{\bf REMARK 6.5}: Lanczos has been trying in vain to do for the $Bianchi$ operator what he did for the $Riemann$ operator, a useless  but possible 
SHIFT BY ONE STEP and to do for the {\it Weyl} operator what he did for the {\it Riemann} operator. However, we shall discover that the dimension $n=4$, which is particularly "{\it fine} " for the classical Killing sequence, is particularly "{\it bad} " for the conformal Killing sequence, a result not known after one century because it cannot be understood without using the {\it Spencer $\delta$-comology} in the following commutative diagram which is explaining therefore what we shall call the "{\it LANCZOS SECRET} ". This diagram allows to consruct the {\it Bianchi} operator ${\cal{D}}_2:F_1 \rightarrow F_2 $ as generating CC for the {\it Riemann} operator ${\cal{D}}_1: F_0=S_2T^* \rightarrow F_1= H^2(g_1)$ defined by a similar diagram and thus only depends on the symbol $g_1$.  \\

\noindent
           \[ \small { \begin{array}{rcccccccccl}
   & 0 &  & 0  & &  0  &  & 0  &  &     \\
   & \downarrow &  &  \downarrow & & \downarrow & & \downarrow & &  & \\
0 \rightarrow & g_4 & \rightarrow &S_4T^*\otimes T& \rightarrow &S_3T^ *\otimes F_0 &\rightarrow & T^*\otimes F_1&\rightarrow & F_2 & \rightarrow 0 \\
 & \downarrow  &  &  \downarrow & & \downarrow & & \parallel & &  \\
0 \rightarrow & T^*\otimes g_3 & \rightarrow &T^*\otimes S_3T^*\otimes T& \rightarrow &T^*\otimes S_2T^ *\otimes F_0 &\rightarrow & T^*\otimes F_1& \rightarrow & 0 &  \\
 & \downarrow  &  &  \downarrow & & \downarrow & & \downarrow & & &  \\
0 \rightarrow &{\wedge}^2 T^*\otimes g_2 & \rightarrow &{\wedge}^2T^*\otimes S_2T^*\otimes T& \rightarrow &{\wedge}^2T^*\otimes T^ *\otimes F_0 &\rightarrow &
 0&&& \\
& \downarrow  &  &  \downarrow & & \downarrow & &  & & & \\
0 \rightarrow &{\wedge}^3 T^*\otimes g_1 & \rightarrow &\underline{{\wedge}^3T^*\otimes T^*\otimes T}& \rightarrow &{\wedge}^3T^*\otimes F_0 &\rightarrow & 0 & &&  \\
  & \downarrow  &  &  \downarrow & & \downarrow & & & & &  \\
0 \rightarrow &{\wedge}^4 T^*\otimes T & = &{\wedge}^4T^*\otimes  T& \rightarrow &0 & & & &&  \\
 & \downarrow  &  &  \downarrow & &  & & & & & \\
 & 0  &  & 0  & &    &  &   &   &  & 
\end{array} } \]   \\  

\noindent
All the vertical down arrows are $\delta$-maps of Spencer and all the vertical columns are exact but the first, which may not be exact only at ${\wedge}^3T^*\otimes g_1$ with cohomology equal to $H^3(g_1)$ because we have:  \\
\[  g_1\simeq {\wedge}^2T^*\subset T^*\otimes T, \,\, g_2=0 \Rightarrow g_3=0 \Rightarrow  g_4=0  \]   

\noindent
A {\it snake-type chase} provides the identification:  \,\,\,  $F_2=H^3(g_1)$  \,\,\,  \fbox {\it SPENCER COHOMOLOGY)}  \\

\noindent
As $g_2=0$, the vector bundle $F_2$ providing the Bianchi identities is defined by the short exact sequence:  \\
\[  \begin{array}{rcccccl}
  0 \longrightarrow & F_2& \longrightarrow  & {\wedge}^3T^*\otimes g_1 & \stackrel{\delta}{\longrightarrow} & {\wedge}^4T^*\otimes T  &\longrightarrow 0  \\
  0 \longrightarrow &20 &\longrightarrow  &24 & \stackrel{\delta}{\longrightarrow} & 4 & \longrightarrow 0  \
\end{array}   \]  \\

When $n=4$, using the duality with respect to the volume form $dx^1\wedge dx^2\wedge dx^3 \wedge dx^4$ in order to change the indices, we obtain successively (care to the signs):  \\
\[  \begin{array}{rccccccccl}
   & B^i_{1,234} & - & B^i_{2,341} & + & B^i_{3,412} & - & B^i_{4,123} & = & 0  \\
   & B_{i1,1}  & - & B_{i2,2} & + & B_{i3,3} & - & B_{i4,4}&  =  &0   \\
i=4 \Rightarrow    &  B_{41,1}& - & B_{42,2}&+ &  B_{43,3} & = & 0 &   \\
  & L_{23,1}& + & L_{31,2}&  + &  L_{12,3} & = &  0
\end{array}  \]   \\  

\noindent
and finally exhibit the {\it Lanczos potential} $L\in {\wedge}^2T^*\otimes T^*$ as a $3$-tensor satisfying: \\

\noindent
\hspace*{3cm}  \fbox{$L_{ij,k} + L_{ji,k}=0, \,\,\,\,\, L_{ij,k} + L_{jk,i} + L_{ki,j}=0 $} \,\,\, $(24-4=20)$  \\  \\

\noindent
{\bf 7) CLASSICAL VERSUS CONFORMAL}  \\  

We obtain successively the following differential sequences for various dimensions:  \\\\
\noindent
CLASSICAL KILLING OPERATOR:   \hspace{3cm}  ${\cal{L}}(\xi)\omega=0 $  \\

\[ \begin{array}{rcccccccccl}
n=2 \,\,\,\,\, \,\, & 2 & \underset 1{\longrightarrow} & 3 & \underset 2 {\longrightarrow}& 1 & \longrightarrow & 0 &  &  &   \\
n=3 \,\,\,\,\,\,\,  & 3&   \underset 1{\longrightarrow} & 6& \underset 2 {\longrightarrow}& 6 & \underset 1{\longrightarrow} & 3 &\longrightarrow  
&  0 &   \\
n=4 \,\,\,\,\,\,\,  & 4 & \underset 1 {\stackrel{K}{\longrightarrow}} &10 & \underset 2{ \stackrel{R}{\longrightarrow}} & 20 & \underset 1 {\stackrel{B}{\longrightarrow}} & 20 & \underset 1 {\longrightarrow}& 6 & \longrightarrow 0 
\end{array}   \]  \\

Euler-Poincar\'{e} characteristic:  $4-10+20-20+6=0$  \\  \\

\noindent
CONFORMAL KILLING OPERATOR:  \hspace{25mm}${\cal{L}}(\xi)\omega=A(x)\omega $\\
\[ metric \,\, density \,\, \Leftrightarrow \,\, {\hat{\omega}}_{ij}={\omega}_{ij} {\mid det(\omega)\mid}^{ - \frac{1}{n}}, \,\, \,\,{\cal{L}}(\xi)\hat{\omega}=0 \]

\[ \begin{array}{rcccccccccccl}
n=3 \,\,\,\,\,\,\,  & 3 & \underset 1{\longrightarrow} & 5 & \underset 3 {\stackrel{?}{\longrightarrow}}& 5 & \underset 1 {\longrightarrow} & 3 & \longrightarrow & 0 & &  &  \\
n=4 \,\,\,\,\, \,\, & 4&   \underset 1{\longrightarrow} & 9 & \underset 2 {\longrightarrow}& 10 & \underset 2{\longrightarrow} & 9 &\underset 1{\longrightarrow } &  4 &\longrightarrow  &0  & \\
n=5 \,\,\,\,\, \,\, & 5 & \underset 1 {\stackrel{CK}{\longrightarrow}} &14 & \underset 2{ \stackrel{W}{\longrightarrow}} & 35 & \underset 1 {\stackrel{?}{\longrightarrow}} & 35 & \underset 2 {\longrightarrow}& 14 & \underset 1 {\longrightarrow} &5 &\longrightarrow 0
\end{array}   \]  \\

Euler-Poincar\'{e} characteristic:  $5-14+35-35+14-5=0$  \\  

We notice that the changes of the successive orders is totally unusual and refer to ([21]) for more details on the computer algebra methods. In particular, when $n=4$, the analogue of the {\it Bianchi} operator is now of order $2$, a result explaining why Lanzos and followers never succeeded adapting the Lanczos tensor potential $L$ for the {\it Weyl} operator. In particular, thanks to Theorems 6.2,6.3 and Corollary 6.4, we have thus solved the {\it Riemann-Lanczos} and {\it Weyl-Lanczos} parametrization problems in arbitrary dimension.\\

\noindent
{\bf 8) VESSIOT STRUCTURE CONSTANTS}  \\

We shall describe the {\it Vessiot structure equations} in the following systematic procedure that can be applied to an arbitrary Lie pseudogroup and ask the reader to compare it with the one adopted in the Introduction. Surprisingly, the example used below has been first exhibited by Vessiot in $1903$ ([29]) and one may refer to ([21,23,24]) for more examples, in particular the case of the Lie pseudogroup of contact transformations when $n=2p+1=3$. \\

\noindent
\fbox{1} \,\, {\bf LIE PEUDOGROUP}: 
\[ y=f(x)\in aut({\mathbb{R}}^2), \Delta (x)=det ({\partial}_if^k(x))\neq 0 \]
\[ \Gamma = \{ y^1= f(x^1), y^2= x^2 / {\partial}_1 f(x^1) \}  \]

\noindent
\fbox{2} /,/,{\bf SYSTEM}:  \,\,\,$y^2dy^1=x^2dx^1 \Rightarrow dy^1\wedge dy^2 = dx^1\wedge dx^2 \Leftrightarrow \Delta = 1 $  \\

\noindent
\fbox{3}\,\, {\bf GENERAL OBJECT}:  \,\,\,  $ \omega =(\alpha,\beta) \in  {\cal{F}} = {\wedge}^1T^* {\times}_X {\wedge}^2T^*  $\\

\noindent
\fbox{4}\,\, {\bf LIE OPERATOR}:  \,\, $ {\cal{D}}\xi\equiv {\cal{L}}(\xi)\omega=0 \Leftrightarrow \{ {\cal{L}}(\xi)\alpha=0, {\cal{L}}(\xi)\beta=0\} $\\

\noindent
\fbox{5}\,\, {\bf MEDOLAGHI EQUATIONS}: 
\[  \{ {\alpha}_r{\partial}_i{\xi}^r+{\xi}^r{\partial}_r{\alpha}_i = 0,  \,\,\, \beta {\partial}_r{\xi}^r+ {\xi}^r{\partial}_r \beta =0 \}\]

\noindent
\fbox{6} \,\, {\bf VESSIOT STRUCTURE EQUATIONS}:\,\,\, \fbox { $ d\alpha=c \beta, \,\,\, c=cst $}\\

\noindent
\fbox{7} \,\, {\bf SPECIAL OBJECT}: \,, $ \alpha=x^2dx^1, \beta = dx^1 \wedge dx^2  \Rightarrow \omega =(x^2,0,1) \Rightarrow c=1$  \\
\[   \bar{\alpha}=dx^1, \bar{\beta} = dx^1 \wedge dx^2  \Rightarrow \bar{\omega} =(1,0,1) \Rightarrow c=0 \Rightarrow \bar{\Gamma}=\{y^1=x^1 + a, y^2=x^2 + f((x^1)\} \] \\

\noindent
\fbox{8}\,\, {\bf DIFFERENTIAL SEQUENCE}:   \hspace{1cm}  $ 0 \rightarrow \Theta \rightarrow \xi \stackrel{{\cal{D}}}{\longrightarrow} \eta \stackrel{{\cal{D}}_1}{\longrightarrow} \zeta  \rightarrow 0 $  \\

\noindent
\[ \left\{ \begin{array}{c} {\xi}^1\\{\xi}^2 \end{array} \right. \longrightarrow \left \{ \begin{array}{rcl}{\alpha}_r{\partial}_i{\xi}^r+{\xi}^r{\partial}_r{\alpha}_i & = & {\eta}^i \\ \beta {\partial}_r{\xi}^r+ {\xi}^r{\partial}_r \beta & = &{\eta}^3 \end{array} \right. \longrightarrow {\partial}_1{\eta}^2 - {\partial}_2{\eta}^1 - c {\eta}^3=\zeta  \]  

\noindent
or, with $D=\mathbb{Q}(x^1,x^2)[d_1,d_2] $ and the rules of differential module theory ([7,15,24]):
\[   0 \rightarrow D \rightarrow D^3 \rightarrow D^2 \rightarrow M \rightarrow 0  \]

\noindent
\fbox{9}\,\, {\bf ADJOINT SEQUENCE}:\hspace{1cm}  $\nu \stackrel{ad({\cal{D}})}{\longleftarrow} \mu \stackrel{ad({\cal{D}}_1)}{\longleftarrow} \lambda $ \\
\[  \lambda  \mid {\partial}_1{\eta}^2 - {\partial}_2{\eta}^1 - c {\eta}^3  \]  

 \[  ad({\cal{D}}_1) \left \{ \begin{array}{rccccl} {\eta}^1 & \rightarrow & {\partial}_2 \lambda & = & {\mu}^1 \\ {\eta}^2 & \rightarrow & - {\partial}_1\lambda & = & {\mu}^2 \\ {\eta}^3 & \rightarrow & - c \lambda & = & {\mu}^3  \end{array} \right.    \hspace{1cm} \fbox {INJECTIVE $\Leftrightarrow  c\neq 0 \Leftrightarrow ext^2(M)=0$}  \]

\[  ad({\cal{D}})  \left \{ \begin{array}{rccccl}  {\xi}^1 & \rightarrow & - {\alpha}_1{\partial}_r{\mu}^r + \beta (c{\mu}^2 - {\partial}_1{\mu}^3) & = & {\nu}^1 \\  {\xi}^2 & \rightarrow & - {\alpha}_2 {\partial}_r {\mu}^r + \beta (-c{\mu}^1- {\partial}_2{\mu}^3) & = & {\nu}^2  \end{array}  \right.\]

\[ \left \{ \begin{array}{lcl}  \bullet\,\,\, c=0 & \Rightarrow & {\partial}_1{\mu}^1 + {\partial}_2{\mu}^2=0, {\mu}^3=0 \\  \bullet \,\,\,  c\neq 0 & \Rightarrow  & {\partial}_1{\mu}^1 + {\partial}_2{\mu}^2=0, {\partial}_1{\mu}^3 - c {\mu}^2=0, {\partial}_2{\mu}^3 + c {\mu}^1=0 \end{array} \right.  \]   \\   

\noindent
In both cases we have $ext^1(M)\neq 0$.  \\  \\

\noindent
{\bf 9) CONCLUSION}  \\

We hope to have convinced the reader that:  \\

\hspace*{2cm}COMPUTER ALGEBRA IS DESPERATELY WANTED \\ 

for studying all these new topics by means of new packages.  \\

\vspace{3cm}

\noindent
{\bf REFERENCES}\\

\noindent
[1]  P. DOLAN, A. GERBER: Janet-Riquier Theory and the Riemann-Lanczos Problem in 2 and 3 Dimensions (2002) arXiv:gr-gq/0212055.  \\
\noindent
[2] S.B. EDGAR: Nonexistence of the Lanczos Potential for the Riemann Tensor in Higher Dimension, Gen. Relat. and Gravitation, 26 (1994) 329-332. \\
\noindent
[3] S.B. EDGAR, A. H\"{O}GLUND: The Lanczos potential for Weyl-Candidate Tensors Exists only in Four Dimension, General Relativity and Gravitation, 32, 12 (2000) 2307. \\
http://rspa.royalsocietypublishing.org/content/royprsa\\
/453/1959/835.full.pdf   \\ 
\noindent
[4] S.B. EDGAR, J.M.M. SENOVILLA: A Local Potential for the Weyl tensor in all dimensions, Classical and Quantum Gravity, 21 (2004) L133.\\
 http://arxiv.org/abs/gr-qc/0408071     \\
\noindent
[5] L.P. EISENHART: Riemannian Geometry, Princeton University Press (1926).  \\
\noindent
[6] M. JANET: Sur les Syst\`{e}mes aux D\'{e}riv\'{e}es Partielles, Journal de Math., 8(3) (1920) 65-151.  \\
\noindent
[7] M. KASHIWARA: Algebraic Study of Systems of Partial Differential Equations, M\'emoires de la Soci\'et\'e Math\'ematique de France 63 
(1995) (Transl. from Japanese of his 1970 Master's Thesis).\\
\noindent
[8] C. LANCZOS: Lagrange Multiplier and Riemannian Spaces, Reviews of Modern Physics, 21 (1949) 497-502.  \\ 
\noindent
[9] C. LANCZOS: The Splitting of the Riemann Tensor, Rev. Mod. Phys. 34 (1962) 379-389.  \\
\noindent
[10] E. MASSA, E. PAGANI: Is the Rieman Tensor Derivable from a Tensor Potential, Gen.Rel. Grav., 16 (1984) 805-816.  \\ 
\noindent
[11] D.G. NORTHCOTT: An Introduction to Homological Algebra, Cambridge University Press (1966).  \\
\noindent
[12] P. O'DONNELL, H. PYE: A Brief Historical Review of the Important Developments in Lanczos Potential Theory, EJTP, 24 (2010) 327-350.  \\
\noindent
[13] J.-F. POMMARET: Systems of Partial Differential Equations and Lie Pseudogroups, Gordon and Breach, New York (1978) 
(Russian translation by MIR, Moscow, 1983) \\
\noindent
[14] J.-F. POMMARET: Dualit\'{e} Diff\'{e}rentielle et Applications, C. R. Acad. Sci. Paris, 320, S\'{e}rie I (1995) 1225-1230.  \\
\noindent
[15] J.-F. POMMARET: Partial Differential Control Theory, Kluwer (2001) 957 pp.\\
\noindent
[16] J.-F. POMMARET: Algebraic Analysis of Control Systems Defined by Partial Differential Equations, in Advanced Topics in Control Systems Theory, Lecture Notes in Control and Information Sciences 311, Chapter 5, Springer (2005) 155-223.\\
\noindent
[17] J.-F. POMMARET: Parametrization of Cosserat Equations, Acta Mechanica, 215 (2010) 43-55.\\
\noindent
[18] J.-F. POMMARET: The Mathematical Foundations of General Relativity Revisited, Journal of Modern Physics, 4 (2013) 223-239.\\
http://dx.doi.org/10.4236/jmp.2013.48A022   \\
\noindent
[19] J.-F. POMMARET: The Mathematical Foundations of Gauge Theory Revisited, Journal of Modern Physics, 5 (2014) 157-170.  \\
http://dx.doi.org/10.4236/jmp.2014.55026    \\
\noindent
[20] J.-F. POMMARET: Airy, Beltrami, Maxwell, Einstein and Lanczos Potentials Revisited, Journal of Modern Physics, 7 (2016) 699-728.  \\
http://dx.doi.org/10.4236/jmp.2016.77068   \\
\noindent
[21] J.-F. POMMARET: Deformation Theory of Algebraic and Geometric Structures, Lambert Academic Publisher, (LAP), Saarbrucken, Germany (2016).  \\
http://arxiv.org/abs/1207.1964  \\
\noindent
[22] J.-F. POMMARET: Why Gravitational Waves Cannot Exist, Journal of Modern Physics, 8,13 (2017) 2122-2158.  \\
http://dx.doi.org/10.4236/jmp.2017.813130   \\
\noindent
[23] J.-F. POMMARET: Homological Solution of the Riemann/Lanczos and Weyl/Lanczos Problems in Arbirary Dimension (2018)  \\
http://arxiv.org/abs/1803.09610    \\
\noindent
[24] J.-F. POMMARET: New Mathematical Methods for Physics, NOVA Science Publisher, New York (2018).  \\
\noindent  
[25] J.-F. POMMARET, A. QUADRAT: Localization and parametrization of linear multidimensional control systems, Systems \& Control Letters, 
37 (1999,)247-260.  \\
\noindent
[26] A. QUADRAT, R. ROBERTZ: A Constructive Study of the Module Structure of Rings of Partial Differential Operators, Acta Applicandae Mathematicae, 133, 187-234 (2014) 187-234. \\
\noindent
[27] J.J. ROTMAN: An Introduction to Homological Algebra, Pure and Applied Mathematics, Academic Press (1979).  \\
\noindent
[28] D.C. SPENCER: Overdetermined Systems of Partial Differential Equations, Bull. Amer. Math. Soc., 75 (1965) 1-114.\\
\noindent
[29] E. VESSIOT: Sur la Th\'{e}orie des Groupes Infinis, Ann. Ec. Normale Sup., 20 (1903) 411-451 (Can be obtained from  http://numdam.org).  \\
\noindent
[30] E. ZERZ: Topics in Multidimensional Linear Systems Theory, Lecture Notes in Control and Information Sciences 256, Springer (2000).\\

\end{document}